\documentclass{amsart}
\usepackage{amsmath, amssymb, graphicx, url}
\usepackage[all]{xy}
\usepackage{type1cm}
\newtheorem{defi}{Definition}[section]
\newtheorem{thm}[defi]{Theorem}

\newtheorem{lem}[defi]{Lemma}
\newtheorem{cor}[defi]{Corollary}

\newtheorem{rem}[defi]{Remark}
\newtheorem{claim}[defi]{Claim}
\newtheorem{conj}[defi]{Conjecture}
\newtheorem{Q}[defi]{Question}
\newtheorem{pf}{proof}

\DeclareMathOperator{\CM}{\underline{CM}}
\DeclareMathOperator{\coh}{\mathrm{coh}}
\DeclareMathOperator{\D}{D^{\mathrm{b}}}
\DeclareMathOperator{\Dperf}{D^{\mathrm{perf}}}
\DeclareMathOperator{\Dsg}{D_{\mathrm{sg}}}
\DeclareMathOperator{\Ext}{\mathrm{Ext}}
\DeclareMathOperator{\gt}{\mathfrak{gt}}

\newcommand{\id}{\mathrm{id}}
\DeclareMathOperator{\Ker}{Ker}
\DeclareMathOperator{\lExt}{\mathcal{E}\mathit{xt}}
\newcommand{\LF}{\mathcal{E}}

\newcommand{\sh}{\mathcal{F}}
\DeclareMathOperator{\Spec}{Spec}
\DeclareMathOperator{\Supp}{Supp}
\newcommand{\stsh}{\mathcal{O}}

\title[dimension of the singularity category]{dimension of the singularity category of a variety with rational singularities}
\author{wahei hara}
\address{Department of Mathematics, School of Science and Engineering, Waseda University, 3-4-1 Ohkubo, Shinjuku, Tokyo 169-8555, Japan}
\email{waheyhey@ruri.waseda.jp}
\date{}

\begin{document}

\begin{abstract}
In this paper, we study the Rouquier dimension of the singularity category of a variety with rational singularity. We construct an upper bound for the dimension of $\Dsg(X)$ if $X$ has at worst rational singularities and $\dim X_{\mathrm{sing}} \leq 1$.
\end{abstract}

\maketitle

\section{Introduction}
The derived category of coherent sheaves on an algebraic variety is an interesting invariant of algebraic varieties and contains much information about the geometry, for example, about its birational properties. Furthermore, We expect that two varieties that share essential features (for example, two varieties that are connected by flops) are expected to have equivalent derived categories.  
If the derived category extracts ``the essential part" of a class of varieties,
one is tempted to work directly with the derived category rather
than individual varieties. When one regards the derived category
as an replacement of a variety, it is natural and important to develop
a dimension theory of triangulated categories.

A notion of a dimension of triangulated categories, called Rouquier dimension, was first introduced by Rouquier in \cite{R08}. We mean the Rouquier dimension of the derived category of an algebraic variety by the Rouquier dimension of the variety for simplicity. Orlov first showed that the Rouquier dimension of a smooth quasi-projective curve over a field is equal to one, and conjectured that the Rouquier dimension of a smooth quasi-projective variety coincides with its Krull dimension (= Conjecture \ref{Orlov conjecture}). Ballard and Favero proved that Orlov's conjecture is true for smooth projective varieties whose derived category has a tilting object that satisfies certain cohomological properties. 

Rouquier showed that the Rouquier dimension of a smooth quasi-projective variety $X$ is not greater than twice the dimension of $X$. 
However, for a variety with singular points, such an upper bound was missing and we only know that the Rouquier dimension is finite.
When we study  the derived category of  a singular variety, one is naturally led to study its singularity categories. The singularity category $\Dsg(X)$ of a singular variety $X$ is a quotient category of the derived category $\D(X)$, and reflects properties of singularities of $X$. Since the theory of the singularity category $\Dsg(X)$ of an affine variety $X = \Spec R$ is related to the theory of Cohen-Macaulay modules over $R$ \cite{B86}, the Rouquier dimension of the singularity category $\Dsg(X)$ is studied by a number of researchers in commutative algebra \cite{AAITY14, AT15, BIKO10, DT15}.

The aim of this paper is to give  an upper bound for $\dim \Dsg(X)$ when $X$ has (at worst) rational singularity, by compering the singularity category $\Dsg(X)$ with the derived category of a resolution $Y$ of $X$ . Specifically, we will prove the following results.

\begin{thm}[= \ref{Main theorem}]
Let $X$ be an algebraic variety over a field $\Bbbk$ such that its category of coherent sheaves has enough locally free. Let $n := \dim X$ and $s := \dim X_{\mathrm{sing}}$, and let us assume
\begin{enumerate}
\item[(a)] there exists a rational resolution $f : Y \to X$, and
\item[(b)] the inequality $\dim f^{-1}(x) \leq n - s$ holds for all $x \in X$.
\end{enumerate}
Then, we have
\begin{enumerate}
\item[(i)] The natural functor $\overline{Rf_*} : \D(Y) \to \Dsg(X)$ has a dense image.
\item[(ii)] If $G \in \D(Y)$ is a generator of $\D(Y)$, then $Rf_*(G)$ is again a generator of $\D(X)$.
\end{enumerate}
\end{thm}

As a corollary, we get an upper bound of $\dim \Dsg(X)$ when $X$ has (at worst) rational singularities with $\dim X_{\mathrm{sing}} \leq 1$.

\begin{cor}[= \ref{cor a}] \label{intro cor2}
Let $X$ be a quasi-projective variety with (at worst) rational singularities of $\dim X_{\mathrm{sing}} \leq 1$, and $Y$ a resolution of $X$. Then, we have an inequality $\dim \Dsg(X) \leq \dim \D(Y) \leq 2\dim X$. 

Moreover, if Conjecture \ref{Orlov conjecture} holds for $Y$, we have $\dim \Dsg(X) \leq \dim X$.
\end{cor}

\vspace{0.1in}

\noindent
\textbf{Acknowledgements. }
The author would like to express his gratitude to his supervisor Professor Yasunari Nagai for beneficial conversations and useful comments. He is also grateful to Professor Shunsuke Takagi, Professor Ryo Takahashi, Doctor Yuki Hirano, Doctor Genki Ouchi for many useful discussions and much helpful advice. 

\section{Definitions and known results}

Let $\mathcal{D}$ be a $\Bbbk$-linear triangulated category.

\begin{defi} \rm
\begin{enumerate}
\item[(1)] Let $\mathcal{I}$ be a full subcategory of $\mathcal{D}$ (not necessarily triangulated). We denote by $\langle \mathcal{I} \rangle$ the smallest full subcategory of $\mathcal{D}$ that contains $\mathcal{I}$ and is closed under shifts, finite direct sums and direct summands. If $\mathcal{I}$ is consisting of only one object $E$, we also write $\langle E \rangle$ instead of $\langle \mathcal{I} \rangle$.
\item[(2)] For two full subcategories $\mathcal{I}_1, \mathcal{I}_2$ of $\mathcal{D}$, we denote by $\mathcal{I}_1 \ast \mathcal{I}_2$ the full subcategory consisting of all objects $M$ that lies in a distinguished triangle $M_1 \to M \to M_2 \to M_1[1]$ where $M_i \in \mathcal{I}_i$.
\item[(3)] For two full subcategories $\mathcal{I}_1, \mathcal{I}_2$ of $\mathcal{D}$, we define $\mathcal{I}_1 \diamond \mathcal{I}_2 := \langle \mathcal{I}_1 \ast \mathcal{I}_2 \rangle$, and we define inductively $\langle \mathcal{I} \rangle_0 := \langle \mathcal{I} \rangle$ and $\langle \mathcal{I} \rangle_k := \langle \mathcal{I} \rangle_{k-1} \diamond \langle \mathcal{I} \rangle_0$. 
\end{enumerate}
\end{defi}

The full subcategory $\langle \mathcal{I} \rangle_k$ was first introduced by Bondal and van den Bergh in \cite{BV03}. Using this notion, Rouquier defined in \cite{R08} a dimension of triangulated categories.

\begin{defi} \rm
Let $G$ be an object of $\mathcal{D}$. We say that $G$ is a \textit{classical generator} of $\mathcal{D}$ if $\bigcup_{i=0}^{\infty} \langle G \rangle_i = \mathcal{D}$. Moreover, We say that $G$ is a \textit{strong generator} of $\mathcal{D}$ if there exists a non-negative integer $g$ such that $\langle G \rangle_g = \mathcal{D}$.
\end{defi}

\begin{defi} \rm
The \textit{generation time} of a strong generator $G$ of $\mathcal{D}$, denoted by $\gt_{\mathcal{D}}(G)$, is the smallest non-negative integer $g$ such that $\langle G \rangle_g = \mathcal{D}$.
\end{defi}

\begin{defi} \rm
The \textit{dimension} of a triangulated category $\mathcal{D}$, denoted by $\dim \mathcal{D}$, is the infimum of the integers $d$ such that there exists a strong generator $G \in \mathcal{D}$ with $d = \gt_{\mathcal{D}}(G)$. If $\mathcal{D}$ have no strong generator, then we set $\dim \mathcal{D} = \infty$.
\end{defi}

The dimension of a triangulated category is sometimes called the \textit{Rouquier dimension} or the \textit{derived dimension}.

\begin{rem} \rm
One can easily show that $\langle G \rangle_i \diamond \langle G \rangle_j = \langle G \rangle_{i+j+1}$. By using this fact, it is immediate to see that if $\mathcal{D}$ has a strong generator then all classical generator is a strong generator.
\end{rem}

\begin{defi} \rm
Let $\Phi : \mathcal{D} \to \mathcal{D}'$ be an exact functor between triangulated categories $\mathcal{D}$ and  $\mathcal{D}'$. We say that $\Phi$ has \textit{dense image} if for all $F \in \mathcal{D}'$ there exists $E \in \mathcal{D}$ such that $F$ is a direct summand of $\Phi(E)$.
\end{defi}

If there exists a functor between triangulated categories $\mathcal{D}$ and $\mathcal{D}'$ that has \textit{dense image}, then we can compare the dimensions of $\mathcal{D}$ and $\mathcal{D}'$.

\begin{lem} \label{dense}
Let $\Phi : \mathcal{D} \to \mathcal{D}'$ be an exact functor between triangulated categories that has dense image. Then, we have $\dim \mathcal{D} \geq \dim \mathcal{D}'$.
\end{lem}

\begin{pf} \rm
Let $G$ be a generator of $\mathcal{D}$ with $\gt_{\mathcal{D}}(G) = \dim \mathcal{D}$, and $F$ an object of $\mathcal{D}'$. Then, there exists an object $E$ of $\mathcal{D}$ such that $F$ is a direct summand of $\Phi(E)$. Then, $\Phi(E) \in \langle \Phi(G) \rangle_{\dim \mathcal{D}}$ and hence $F \in \langle \Phi(G) \rangle_{\dim \mathcal{D}}$. \qed
\end{pf}

It is a natural and fundamental question to calculate the dimension of bounded derived categories of algebraic varieties. The following conjecture is due to Orlov in \cite{O09b}.

\begin{conj}[Orlov \cite{O09b}] \label{Orlov conjecture}
Let $X$ be a smooth quasi-projective variety over a field $\Bbbk$. Then we have
\[ \dim X = \dim \D(X). \]
\end{conj}

In relation to this conjecture, the following results are known:

\begin{thm}[Rouquier \cite{R08}] \label{thm Rouquier} Let $\Bbbk$ be an arbitrary field. Then : 
\begin{enumerate}
\item[(i)] Conjecture \ref{Orlov conjecture} holds for smooth affine schemes of finite type over $\Bbbk$.
\item[(ii)] If $X$ is a reduced separated scheme of finite type over $\Bbbk$, then we have $\dim \D(X) \geq \dim X$.
\item[(iii)] If $X$ is a smooth quasi-projective scheme over $\Bbbk$, then we have $\dim \D(X) \leq 2\dim X$.
\item[(iv)] If $X$ be a separated scheme of finite type over $\Bbbk$, and if $\Bbbk$ is a perfect field, then we have $\dim \D(X) < \infty$.
\item[(v)] Conjecture \ref{Orlov conjecture} is true for flag varieties of type A and smooth quadric hypersurfaces.
\end{enumerate}
\end{thm}

\begin{thm}[Orlov \cite{O09b}]
Conjecture \ref{Orlov conjecture} holds for smooth quasi-projective curves over a filed $\Bbbk$.
\end{thm}

Moreover, Ballard and Favero showed in \cite{BF12} that Conjecture \ref{Orlov conjecture} holds for smooth projective varieties whose derived categories have a tilting object $T$ that satisfies $\Ext^i(T, T \otimes \omega_X^{-1}) = 0$ for all $i > 0$.

The purpose of this paper is to study the dimension of derived categories of singular varieties. For a singular variety $X$, the inequality $\dim \D(X) \leq 2\dim X$ does not hold in general.

When we study the derived category of a singular variety $X$, the \textit{singularity category} $\Dsg(X)$ naturally appears as the difference between the derived category $\D(X)$ and the category of perfect complexes $\Dperf(X)$.

\begin{defi} \rm
 Let $X$ be a noetherian scheme whose category of coherent sheaves has enough locally free sheaves. We say a complex $F$ in $\D(X)$ is perfect if $F$ is locally isomorphic to a bounded complex of locally free sheaves.  We denote by $\Dperf(X)$ the triangulated subcategory of $\D(X)$ consisting of perfect complexes. The \textit{singularity category} of $\D(X)$, denoted by $\Dsg(X)$, is by definition the quotient category
\[ \Dsg(X) := \D(X)/ \Dperf(X). \]
If $X$ is an affine scheme $X = \Spec R$, we also write $\Dsg(R) = \Dsg(X)$.
\end{defi}

For the definition of the quotient category, see \cite{O04}. In the context of ring theory, the singularity category $\Dsg(R)$ is also called the \textit{stable derived category} of $R$. If $R$ is a Gorenstein ring, it is equivalent to the \textit{stable category of (maximal) Cohen-Macaulay modules} $\CM(R)$. Hence the study of the singularity category $\Dsg(X)$ is related to the study of Cohen-Macaulay modules.   For the details, see \cite{B86}.

\begin{rem} \rm
The concept of singularity category was first introduced by Orlov in \cite{O04}. He showed there that if $U$ is an open subscheme of $X$ such that it contains the singular locus $X_{\mathrm{sing}}$ of $X$, then we have a natural equivalence of triangulated categories $\Dsg(U) \simeq \Dsg(X)$. Thus, we can say that the singularity category $\Dsg(X)$ of $X$ is a category that reflects the properties of the singularity of $X$.
\end{rem}

Of course, we have $\dim \Dsg(X) \leq \dim \D(X)$ by definition, and by Theorem \ref{thm Rouquier} (iv), we have $\dim \Dsg(X) < \infty$ if $X$ is a separated algebraic scheme over a perfect field. If $X = \Spec R$ where $R$ is complete Cohen-Macauray local ring with an isolated singularity, Dao and Takahashi constructed an upper bound of $\dim \Dsg(X)$ that depends on some invariants determined by the Jacobian ideal $J$ of $R$ (see \cite{DT15} for details).



Later, we will show in Corollary \ref{cor a} that $\dim \Dsg(X)$ has an upper bound of another type when $X$ has at worst rational singularities and $\dim X_{\mathrm{sing}} \leq 1$. Our bound for $\Dsg(X)$ depends only on the Krull dimension of $X$.

For a lower bound of $\Dsg(X)$, Aihara and Takahashi showed the following  theorem in \cite{AT15}. Recall that a local ring $(R, \mathfrak{m})$ is called \textit{complete intersection} if its $\mathfrak{m}$-adic completion $\hat{R}$ is a residue ring of a complete regular local ring by an ideal generated by a regular sequence. A prime ideal $\mathfrak{p}$ of a ring $R$ is called \textit{complete intersection} if the localization $R_{\mathfrak{p}}$ of $R$ by $\mathfrak{p}$ is complete intersection. The \textit{embedding dimension} of a local ring $(R, \mathfrak{m})$, written by $\mathrm{emb.dim}(R)$, is the dimension of $\mathfrak{m}/\mathfrak{m}^2$ as a vector space over the residue field $R/\mathfrak{m}$. The \textit{codimension} of a local ring $(R, \mathfrak{m})$ is defined by $\mathrm{codim}(R) := \mathrm{emb.dim}(R) - \dim R$.

\begin{thm}[Aihara-Takahashi \cite{AT15}] \label{AT}
Let $R$ be a ring. Then we have the inequality
\[ \dim \Dsg(R) \geq \sup\{ \mathrm{codim}(R_{\mathfrak{p}}) - 1 \mid \mathfrak{p} \in \mathrm{CI}(R) \} \]
where $\mathrm{CI}(R)$ is a complete intersection locus of $R$, that is, the set of all complete intersection prime ideals.
\end{thm}

Later, we will use this result to show an upper bound of the embedding dimension of a singular point that is rational (Remark \ref{rem b}).

\section{Results}

\begin{defi} \rm
Let $X$ be an algebraic variety over a field $\Bbbk$.  We say that a partial resolution $f : Y \to X$ is \textit{rational} if it satisfies $Rf_*\stsh_Y \simeq \stsh_X$.
\end{defi}

Note that we do not assume the Grauert-Riemenschneider vanishing in the definition of  a rational resolution even in the case when $\mathrm{char}(\Bbbk) > 0$. 

The following is the main result in this paper.

\begin{thm} \label{Main theorem}
Let $X$ be an algebraic variety over a field $\Bbbk$ whose category of coherent sheaves has enough locally free. Let $n := \dim X$ and $s := \dim X_{\mathrm{sing}}$, and let us assume that
\begin{enumerate}
\item[(a)] There exists a rational partial resolution $f : Y \to X$, 
\item[(b)] $X$ is Cohen-Macaulay, and
\item[(c)] the inequality $\dim f^{-1}(x) \leq n - s$ holds for all $x \in X$.
\end{enumerate}
Then, we have
\begin{enumerate}
\item[(i)] The natural functor $\overline{Rf_*} : \D(Y) \to \Dsg(X)$ has a dense image.
\item[(ii)] If $G \in \D(Y)$ is a  strong generator of $\D(Y)$, then $Rf_*(G)$ is again a strong generator of $\D(X)$.
\end{enumerate}
\end{thm}

\begin{pf} \rm
Let $\bar{F}$ be a non-zero object in $\Dsg(X)$ and $F$ be a lift of $\bar{F}$ to $\D(X)$. Let $\LF \in D^-(X)$ be a locally free resolution of $F$ in $D^-(X)$. For sufficiently large $k$, we consider a truncation $\sigma^{\geq -k} \LF$ and the natural map $\sigma^{\geq -k} \LF \to F$. Here, the truncation $\sigma^{\geq -k} \LF$ means a complex
\[ \cdots \to 0 \to  0 \to \mathcal{E}^{-k} \xrightarrow{d} \mathcal{E}^{-k+1} \xrightarrow{d} \mathcal{E}^{-k+2} \xrightarrow{d} \cdots. \]
Extending the map $\tau^{\geq -k} \LF \to F$ to the distinguished triangle in $\D(X)$, we have a distinguished triangle
\[ \sh_k[k] \to \sigma^{\geq -k} \LF \to F \to \sh_k[k+1], \]
where $\sh_k$ is a coherent sheaf on $X$. Note that $\sigma^{\geq -k} \LF$ is a perfect complex in $\D(X)$, hence we have $\bar{\sh_k}[k+1] \simeq \bar{F}$ in $\Dsg(X)$. We also note that there exists an exact sequence
\[ 0 \to \sh_k \to \LF^{-k} \to \sh_{k-1} \to 0 \]
for sufficiently large $k$.

Next, we consider the \textit{underived} pullback $f^*\sh_k \in \coh(Y)$ and its derived push-forward $Rf_*(f^*\sh_k)$. We have natural morphisms $\sh_k \to f_*f^*\sh_k$ and $f_*f^*\sh_k \to Rf_*(f^*\sh_k)$, and we have a composition morphism $\sh_k \to Rf_*(f^*\sh_k)$. Extending this morphism to the distinguished triangle in $\D(X)$, we have
\[ \sh_k \to Rf_*(f^*\sh_k) \to C_k \xrightarrow{\phi} \sh_k[1]. \]
Now, we need the following claim.

\begin{claim} \label{claim a}
The natural morphism $\sh_k \to f_*f^*\sh_k$ is injective and split.
\end{claim}

\begin{pf} \rm
Let us consider the exact sequence
\[ 0 \to \sh_k \to \LF^{-k} \xrightarrow{\partial^{-k}} \LF^{-k+1}. \]
Let us apply the functor $f_*f^*$ for this exact sequence. From the projection formula, we have $f_*f^*\LF^i \simeq \LF^i \otimes f_*\stsh_Y \simeq \LF^i$. Furthermore, since the functor $f_*$ is left exact, we have an isomorphism $\sh_k \simeq f_*\Ker(f^*\partial^{-k})$. Thus, we have the following commutative diagram
\[ \begin{xy} \xymatrix{
\sh_k \ar[r]^>>>>>>{\text{adj}} \ar@/_3mm/[rrd]_{\id} & f_*f^*\sh_k \ar[r] \ar@{}[rd]|{\circlearrowright} &  f_*\Ker(f^*\partial^{-k}) \ar[r] \ar@{}[rd]|{\circlearrowright} & f_*f^*\LF^{-k} \ar[r]^>>>>{f_*f^*\partial^{-k}} \ar@{}[rd]|{\circlearrowright} & f_*f^*\LF^{-k+1} \\
 & & \sh_k \ar[r] \ar@{}[u]|{\rotatebox{90}{$\simeq$}} & \LF^{-k} \ar[r]_{\partial^{-k}} \ar@{}[u]|{\rotatebox{90}{$\simeq$}} & \LF^{-k+1}. \ar@{}[u]|{\rotatebox{90}{$\simeq$}}
} \end{xy} \]
Hence the morphism $\mathrm{adj} : \sh_k \to f_*f^*\sh_k$ is injective and split. \qed
\end{pf}

By this claim, we have
\begin{align*}
\mathcal{H}^i(C_k) = \begin{cases}
f_*f^*\sh_k/\sh_k &(i=0) \\
R^if_*(f^*\sh_k) &(i \neq 0).
\end{cases}
\end{align*}
Next, we will show that $\Ext_X^1(C_k, \sh_k)$ is zero. To see this, we first prove the following claim.

\begin{claim} \label{claim b}
Let $T$ be a coherent sheaf on $X$ such that $\dim \Supp(T) \leq d$. Then we have
\[ \Ext_X^i(T, \mathcal{P}) = 0 \text{ if $i < n - d$} \]
for any locally free sheaf $\mathcal{P}$ on $X$.
\end{claim}

We postpone the proof of Claim \ref{claim b} until we finish the proof of Theorem \ref{Main theorem}.
Next, we prove the following claim.

\begin{claim} \label{claim c}
Let $T$ be a coherent sheaf on $X$ such that $\dim \Supp(T) \leq d$. Then we have
\[ \Ext_X^i(T, \sh_k) = 0 \text{ for $i < n - d$,} \]
if the integer $k$ is sufficiently large.
\end{claim}

\begin{pf} \rm
Applying the functor $\Ext^i(T, -)$ to the exact sequence
\[ 0 \to \sh_k \to \LF^{-k} \to \sh_{-k+1} \to 0 \]
and using Claim \ref{claim b}, we have an isomorphism
\[ \Ext^i_X(T, \sh_k) \simeq \Ext^{i-1}_X(T, \sh_{k-1}) \text{ for $i < n -d$.} \]
If the integer $k$ is sufficiently large, we can repeatedly use this isomorphism, and we have
\[ \Ext_X^i(T, \sh_k) \simeq \Ext_X^{-1}(T, \sh_{k-i-1}) = 0 \]
for all $i < n -d$. \qed
\end{pf}

Finally, we need the following claim.

\begin{claim} \label{claim d}
Let $\sh$ be a coherent sheaf on $X$. Then $R^if_*(f^*\sh) = 0$ for all $i > n - s - 2, i \neq 0$.
\end{claim}

The proof of Claim \ref{claim d} is also put off until we complete the proof of Theorem \ref{Main theorem}.
Now we consider the spectral sequence
\[ E_2^{p, q} = \Ext_X^p(\mathcal{H}^{-q}(C_k), \sh_k) \Rightarrow \Ext_X^{p+q}(C_k, \sh_k). \]
By Claim \ref{claim c}, we have $\Ext_X^p(\mathcal{H}^{-q}(C_k), \sh_k) = 0$, for all $p < n -s$ and for sufficiently large $k$, because $\Supp(\mathcal{H}^i(C_k))$ is contained in $X_{\mathrm{sing}}$ for all $i \in \mathbb{Z}$. By Claim \ref{claim d}, we have $\Ext_X^p(\mathcal{H}^{-q}(C_k), \sh_k) = 0$ for all $q < - n + s + 2$. Thus, we have
\[ \Ext_X^p(\mathcal{H}^{-q}(C_k), \sh_k) = 0 \text{ for $p+q \leq 1$,} \]
and from this, we have $\Ext_X^1(C_k, \sh_k) = 0$ for sufficiently large $k$.
Therefore, the distinguished triangle
\[ \sh_k \to Rf_*(f^*\sh_k) \to C_k \to \sh_k[1] \]
is split, and we have $Rf_*(f^*\sh_k) \simeq \sh_k \oplus C_k$. This shows (i).

Let $G \in \D(Y)$ be a generator of $\D(Y)$ such that $\gt_{\D(Y)}(G) = g$. Then we have $f^*\sh_k[k], f^*\tau^{\geq -k}\LF \in \langle G \rangle_g$. Applying the functor $Rf_*$, we have $Rf_*(f^*\sh_k)[k] \in \langle Rf_*G \rangle_g$ and $Rf_*(f^*\sigma^{\geq -k}\LF) \simeq \sigma^{\geq -k}\LF \in \langle Rf_*G \rangle_g$. Because $\sh_k[k]$ is a direct summand of $Rf_*(f^*\sh_k)[k]$, we also have $\sh_k[k] \in \langle Rf_*G \rangle_g$. From the distinguished triangle
\[ \sh_k[k] \to \sigma^{\geq -k}\LF \to F \to \sh_k[k+1], \]
we have
\[ F \in \langle Rf_*G \rangle_g \diamond \langle Rf_*G \rangle_g = \langle Rf_*G \rangle_{2g+1}. \]
This proves (ii) and the proof of Theorem \ref{Main theorem} was completed. \qed
\end{pf}

Let us prove the remaining claims.

\begin{pf}[of Claim \ref{claim b}] \rm
Let us consider the local extension shaves $\lExt_X^i(T, \mathcal{P})$ for each $i \geq 0$. Let $x \in \Supp(T)$ be a point and $(R, \mathfrak{m}) := (\stsh_{X, x}, \mathfrak{m}_x)$ the local ring of $X$ at $x$. Then we have
\[ \lExt_X^i(T, \mathcal{P})_x \simeq \Ext_R^i(T_x, \mathcal{P}_x) = \Ext^i(T_x, R)^{\oplus \mathrm{rk}(\mathcal{P})} . \]

We want to show that $\Ext_R^i(T_x, R) = 0$ for $i < n - d$. First, if $T_x = R/\mathfrak{m}$, then this vanishing follows from the definition of the depth. Next, by the induction on the dimension $d$, we will show the vanishing if $T_x = R/\mathfrak{p}$ where $\mathfrak{p}$ is any prime ideal of $R$. Let us consider the following exact sequence
\[ 0 \to R/\mathfrak{p} \xrightarrow{\cdot a} R/\mathfrak{p} \to R/(\mathfrak{p}, a) \to 0 \]
where $a \in \mathfrak{m}$. From the induction hypothesis, we have $\Ext^i(R/(\mathfrak{p}, a), R) = 0$ for $i < n - d - 1$. Thus, we have  $a \cdot \Ext^i(R/\mathfrak{p}, R) = \Ext^i(R/\mathfrak{p}, R)$ for a $x \in \mathfrak{m}$ and for all $i < n - d$. Hence, by Nakayama's lemma, we have $\Ext^i(R/\mathfrak{p}, R) = 0$ for $i < n -d$.
Finally, by taking a filtration of $T_x$
\[ 0 \subset T_m \subset T_{m-1} \subset \cdots \subset T_2 \subset T_1 \subset T_0 = T_x \]
such that $T_j / T_{j+1} \simeq R/\mathfrak{p}_j$ where $\mathfrak{p}_j$ is a prime ideal of $R$, we can check the vanishing
\[ \Ext_R^i(T_x, R) = 0 ~~~ (i < n-d)\]
for general $T_x$ with $\dim \Supp(T_x) \leq d$.

Thus, we have $\lExt_X^i(T, \mathcal{P}) = 0$ for $i < n - d$. Now, we consider the local-to-global spectral sequence
\[ E^{p. q}_2 = H^p(X, \lExt^q_X(T, \mathcal{P}) )\Rightarrow \Ext_X^{p+q}(T, \mathcal{P}). \]
Since $H^p(X, \lExt_X^q(T, \mathcal{P})) = 0$ for $p+q < n - d$, we have $\Ext_X^{p+q}(T, \mathcal{P}) = 0$ for $p+q < n - d$. \qed
\end{pf}

\begin{pf}[of Claim \ref{claim d}] \rm
Since we assumed that $\dim f^{-1}(x) \leq n - s$ for all $x \in X$, we have $R^if_*\mathcal{G} = 0$ for any coherent sheaf $\mathcal{G}$ on $Y$ and for all $i > n - s$. Thus, we only need to check $R^{n-s}f_*(f^*\sh) = 0$ and $R^{n-s-1}f_*(f^*\sh) = 0$. We take a locally free resolution $\LF' \to \LF \to \sh \to 0$. Then we have an exact sequence
\[ 0 \to K \to f^*\LF \to f^*\sh \to 0 \]
and a surjection $f^*\LF' \twoheadrightarrow K$. By the projection formula, we have $R^if_*(f^*\mathcal{P}) \simeq \mathcal{P} \otimes R^if_*\stsh_Y = 0$ for all $i > 0$ and for any locally free sheaf $\mathcal{P}$ on $X$.  Therefore, we have $R^{n -s}f_*(f^*\sh) = 0$ and $R^{n-s}f_*K = 0$. Moreover, from the exact sequence
\[ R^{n-s-1}f_*(f^*\LF) \to R^{n-s-1}f_*(f^*\sh) \to R^{n-s}f_*K, \]
we have $R^{n-s-1}f_*(f^*\sh) = 0$. \qed
\end{pf}

\begin{cor} \label{cor a}
Let $X$ be a quasi-projective variety over a field $\Bbbk$ of characteristic zero. Assume that $X$ has at worst rational singularities of $\dim X_{\mathrm{sing}} \leq 1$. Let $Y$ be a resolution of $X$. Then,  we have $\dim \Dsg(X) \leq \dim \D(Y) \leq 2\dim X$.
Moreover, if Conjecture \ref{Orlov conjecture} holds for $Y$, we have $\dim \Dsg(X) \leq \dim X$.
\end{cor}

\begin{pf} \rm
Note that a rational singularity over a field of characteristic zero is automatically Cohen-Macaulay and hence any resolution $f : Y \to X$ satisfies the assumption in Theorem \ref{Main theorem}. Hence the first inequality $\dim \Dsg(X) \leq \dim \D(Y)$ follows from Theorem \ref{Main theorem} and Lemma \ref{dense}. Since $Y$ is smooth, we have the second inequality $\dim \D(Y) \leq 2 \dim Y = 2 \dim X$ by Theorem \ref{thm Rouquier}.
\qed 
\end{pf}

\begin{rem} \label{rem a} \rm
As a corollary of the proof of Theorem \ref{Main theorem}, we also have an inequality $\dim \D(X) \leq 2\dim \D(Y) + 1$. However, this bound is not optimal. Indeed, Favero showed $\dim \D(X) \leq \dim \D(Y)$ under the assumption that $X$ has at worst rational singularity without an assumption on $\dim X_{\mathrm{sing}}$ and $Y$ is a resolution of $X$, in his PhD thesis \cite{Fa09}. The inequality $\dim \Dsg(X) \leq \dim \D(Y)$ as in Corollary \ref{cor a} also follows from his result, since we trivially have $\dim \Dsg(X) \leq \dim \D(X)$. His approach is based on so-called ``Ghost-Lemma" (see \cite{BFK12}), which is quite different from our argument via the density of the functor $\overline{Rf_*}$ .
\end{rem}

\begin{rem} \label{rem b} \rm
By combining the result in Remark \ref{rem a} with the result of Aihara and Takahashi (Theorem \ref{AT}), we can partly recover a well-known inequality in the singularity theory from the Orlov's conjecture : 

\textit{Let $X$ be a quasi-projective variety with (at worst) rational singularities. 
We assume that Conjecture \ref{Orlov conjecture} holds for some resolution of $X$. Then, we have the inequality
\[ \mathrm{codim}(\stsh_{X, x}) \leq \dim X + 1 \]
for all complete intersection point $x \in X$.}

In the case if $\dim X_{\mathrm{sing}} \leq 1$, the same inequality also follows from our Corollary \ref{cor a}.

It is known that there is a stronger inequality that holds for complete intersection rational singularities over an algebraic closed field of characteristic zero. Namely, the following inequality holds for a local ring  $(R, \mathfrak{m})$ essencially of finite type over $\Bbbk$ that is complete intersection and has a rational singularity: $\mathrm{codim}(R) \leq \dim R - 1$. The author learned this inequality from Professor S. Takagi in our private communication. For the sake of completeness, we give the proof of this inequality here, which is also according to Professor S. Takagi. 

Let $\Bbbk$ be an algebraically closed field of characteristic zero and $X=\mathbb{A}^e_{\Bbbk}=\Spec \Bbbk[X_1, \dots, X_e]$. 
Let $I=(f_1, \dots, f_c), \mathfrak{m}=(X_1, \dots, X_e) \subset \Bbbk[X_1, \dots, X_e]$. Suppose that each $f_i$ is contained in $\mathfrak{m}^2$ and  $Z:=\Spec \Bbbk[X_1, \dots, X_e]/I$ has only rational singularities and is locally complete intersection in $X$.
Let $f:Y \to X$ be the blow-up of $X$ along $Z$ with exceptional divisor $F$ and let $g: X' \to Y$ be a log resolution of $(Y, F, \mathfrak{m}\stsh_Y)$.
We can write $K_{X'/X}=\sum_i a_i E_i+(c-1)F'$ and $g^*F=\sum_{i} b_i E_i+F'$, where the $E_i$ are $g$-exceptional divisors and $F'$ is the strict transform of $F$. By the definition of log resolution, $\mathfrak{m}\stsh_{X'}$ is invertible, so we can write $\mathfrak{m}\stsh_{X'}=\stsh_{X'}(-\sum_i m_i E_i)$.
Since $Z$ has canonical singularities, $a_i \geq c b_i$ for all $i$.
(see for example \cite[Theorem 2.1]{Mu01}).
On the other hand,  since $I$ is contained in $\mathfrak{m}^2$, one has $b_i \geq 2 m_i$ for all $i$. Thus, $a_i \geq 2c m_i$ for all $i$, which implies that $(X, \mathfrak{m}^{2c})$ is klt (we say that $(X, \mathfrak{m}^t)$ is klt if the coefficients of $K_{X'/X}-\sum_i t m_i E_i$ is greater than or equal to $-1$).
Thus, the log canonical threshold $\mathrm{lct}(\mathfrak{m})$ of $\mathfrak{m}$ is greater than $2c$.
On the other hand, it is easy to see that $\mathrm{lct}(\mathfrak{m})=e$. Thus, we have the desired inequality $e \geq 2c+1$.

In this proof, we look at the relative canonical bundle of $f : X' \to X$, while in our proof of Theorem \ref{Main theorem}, we compare the derived category of $\D(X)$ and the derived category $\D(Y)$ of the resolution $Y$. Here appears the ubiquitous analogy between the geometry of canonical divisors and the theory of derived categories that is widely observed (see for example \cite{K09}).
\end{rem}

If we believe the conjectural upper bound $\dim \D(X) \leq \dim X$, it is also important and interesting to ask the following questions:

\begin{Q} \rm
Does there exist a quasi-projective variety $X$ with (at worst) rational singularities that satisfies $\dim \Dsg(X) = \dim X$ ? If not, what is the optimal bound for it? What kind of singularity does $X$ have if $\dim \Dsg(X)$ attains the maximum?
\end{Q}


\begin{thebibliography}{99}

\bibitem[AAITY14]{AAITY14} T. Aihara, T. Ayara, O. Iyama, R. Takahashi, \textit{Dimensions of triangulated categories with respect to subcategories}, J. Algebra \textbf{399} (2014), 205--219.

\bibitem[AT15]{AT15} T. Aihara, R. Takahashi, \textit{Generators and dimensions of derived categories of modules}, Comm. Algebra \textbf{43} (11) (2015), 5003--5029.

\bibitem[BF12]{BF12} M. Ballard, D. Favero, \textit{Hochschild dimensions of tilting objects}, Int. Math. Res. Not. \textbf{11} (2012), 2607--2645.

\bibitem[BFK12]{BFK12} M. Ballard, D. Favero, L. Katzarkov, \textit{Orlov spectra: bounds and gaps}, Invent. Math. \textbf{189}(2) (2012), 359--430.


\bibitem[BIKO10]{BIKO10} P.A. Bergh, S.B. Iyenger, H. Krause, S. Oppermann, \textit{Dimensions of triangulated categories via {K}oszul objects}, Math. Z. \text{265}(4) (2010), 849--864.

\bibitem[BV03]{BV03} A. Bondal, M. van den Bergh, \textit{Generators and representability of functors in commutative and noncommutative geometry}, Mosc. Math. J. \textbf{3} (1) (2003), 1--36, 258.


\bibitem[Bu86]{B86} R.-O. Buchweitz, \textit{Maximal Cohen-Macaulay modules and Tate-cohomology over Gorenstein rings}, unpublished paper, \url{http://hdl.handle.net/1807/16682}.



\bibitem[DT15]{DT15} H. Dao, R. Takahashi, \textit{Upper bounds for dimensions of singularity categories}, C. R. Math. Acad. Sci. Paris, \textbf{353}(4) (2015), 297--301.


\bibitem[Fa09]{Fa09} D. Favero, \textit{A study of the geometry of the derived category}, unpublished PhD thesis, University of Pennsylvania, 2009.

\bibitem[Ka09]{K09} Y. Kawamata, \textit{Derived categories and birational geometry}, Algebraic geometry---{S}eattle 2005. {P}art 2, Proc. Sympos. Pure Math. \textbf{80} (2009), 655--665.


\bibitem[Mu01]{Mu01} M. Musta{\c{t}}{\u{a}}, \textit{Jet schemes of locally complete intersection canonical singularities}, Invent. Math. \textbf{145}(3) (2001), 397--424.

\bibitem[Or04]{O04} D. Orlov, \textit{Triangulated categories of singularities and D-branes in Landau-Ginzburg models}, Proc. Steklov Inst. Math, \textbf{246} (2004), 227--248.

\bibitem[Or09]{O09b} D. Orlov, \textit{Remarks on generators and dimension of triangulated categories}, Mosc. Math. J. \textbf{9} (1) (2009), 153--159.

\bibitem[Ro08]{R08} R. Rouquier, \textit{Dimension of triangulated categories}, J. K-Theory \textbf{1} (2) (2008) 193--256.


\end{thebibliography}
\end{document}